\documentclass[12TP,draft]{article}

\begin{document}

\begin{center}
\LARGE\noindent\textbf{ On longest non-Hamiltonian Cycles in Digraphs with the Conditions of Bang-Jensen, Gutin and Li}\\

\end{center}
\begin{center}
\noindent\textbf{S.Kh. Darbinyan and I.A. Karapetyan}\\

Institute for Informatics and Automation Problems, Armenian National Academy of Sciences

E-mails: samdarbin@ipia.sci.am, isko@ipia.sci.am\\
\end{center}

\textbf{Abstract}\\

 Let $D$ be a strong digraph on $n\geq 4$ vertices. In [2, J. Graph Theory 22 (2) (1996) 181-187)], J. Bang-Jensen, G. Gutin and H. Li proved the following theorems: If (*) $d(x)+d(y)\geq 2n-1$ and $min \{ d(x), d(y)\}\geq n-1$ for every pair of non-adjacent vertices $x, y$ with a common in-neighbour or (**) $min \{ d^+(x)+ d^-(y),d^-(x)+ d^+(y)\}\geq n$ for every pair of non-adjacent vertices $x, y$ with a common in-neighbour or a common out-neighbour, then $D$ is hamiltonian.
 In this paper we show that: (i) if $D$  satisfies the condition (*) and the minimum semi-degree of $D$ at least two or (ii) if $D$ is not directed cycle and satisfies the condition (**), then either $D$ contains a cycle of length $n-1$ or $n$ is even and $D$ is isomorphic to complete bipartite digraph or to complete bipartite digraph minus one arc.  \\


\noindent\textbf{1. Introduction and Terminology}\\

We shall assume that the reader is familiar with the standard terminology on directed graphs (digraphs) and refer the reader to monograph of Bang-Jensen and Gutin [1] for terminology not discussed here. In this paper we consider finite digraphs without loops and multiple arcs. For a digraph $D$, we denote by $V(D)$ the vertex set of $D$ and by  $A(D)$ the set of arcs in $D$. Often we will write $D$ instead of $A(D)$ and $V(D)$. The arc of a digraph $D$ directed from $x$ to $y$ is denoted by $xy$. For disjoint subsets $A$ and  $B$ of $V(D)$  we define $A(A\rightarrow B)$ \, as the set $\{xy\in A(D) / x\in A, y\in B\}$ and $A(A,B)=A(A\rightarrow B)\cup A(B\rightarrow A)$. If $x\in V(D)$ and $A=\{x\}$ we write $x$ instead of $\{x\}$. If $A$ and $B$ are two disjoint subsets of $V(D)$ such that every vertex of $A$ dominates every vertex of $B$, then we say that $A$ dominates $B$, denoted by $A\rightarrow B$. The out-neighbourhood of a vertex $x$ is the set $N^+(x)=\{y\in V(D) / xy\in A(D)\}$ and $N^-(x)=\{y\in V(D) / yx\in A(D)\}$ is the in-neighbourhood of $x$. Similarly, if $A\subseteq V(D)$ then $N^+(x,A)=\{y\in A / xy\in A(D)\}$ and $N^-(x,A)=\{y\in A / yx\in A(D)\}$. We call the vertices in $N^+(x)$, $N^-(x)$, the out-neighbours and in-neighbours of $x$. The out-degree of $x$ is $d^+(x)=|N^+(x)|$ and $d^-(x)=|N^-(x)|$ is the in-degree of $x$. The out-degree and in-degree of $x$ we call its semi-degrees. Similarly, $d^+(x,A)=|N^+(x,A)|$ and $d^-(x,A)=|N^-(x,A)|$. The degree of the vertex $x$ in $D$ defined as $d(x)=d^+(x)+d^-(x)$ (similarly, $d(x,A)=d^+(A)+d^-(A)$). The subdigraph of $D$ induced by a subset $A$ of $V(D)$ is denoted by $\langle A\rangle$. The path (respectively, the cycle) consisting of the distinct vertices $x_1,x_2,\ldots ,x_m$ ( $m\geq 2 $) and the arcs $x_ix_{i+1}$, $i\in [1,m-1]$  (respectively, $x_ix_{i+1}$, $i\in [1,m-1]$, and $x_mx_1$), is denoted  $x_1x_2\cdots x_m$ (respectively, $x_1x_2\cdots x_mx_1$). For a cycle  $C_k=x_1x_2\cdots x_kx_1$, the subscripts considered modulo $k$, i.e. $x_i=x_s$ for every $s$ and $i$ such that  $i\equiv s\, (\hbox {mod} \,k)$. If $P$ is a path containing a subpath from $x$ to $y$ we let $P[x,y]$ denote that subpath. Similarly, if $C$ is a cycle containing vertices $x$ and $y$, $C[x,y]$ denotes the subpath of $C$ from $x$ to $y$.
 A digraph $D$ is strongly connected (or just strong) if there exists a path from $x$ to $y$ and a path from $y$ to $x$ in $D$ for every choice of distinct vertices $x$,\,$y$ of $D$. A digraph $D$ is semicomplete if, for every pair of distinct vertices $x$ and $y$, there is at least one arc between them and is locally semicomplete, if $\langle N^+(x)\rangle$ and $\langle N^-(x)\rangle$ are both semicomplete for every $x$ of $D$. We will denote the  complete bipartite digraph  with partite sets of cardinalities $p$, $q$ by $K^*_{p,q}$. Two distinct vertices $x$ and $y$ are adjacent if $xy\in A(D)$ or $yx\in A(D) $ (or both). We denote by $a(x,y)$ the number of arcs between the vertices $x$ and $y$. In particular, $a(x,y)=0$ (respectively, $a(x,y)\not=0$) means that $x$ and $y$ are not adjacent (respectively, are adjacent). 

For integers $a$ and $b$, $a\leq b$, let $[a,b]$  denote the set of all integers which are not less than $a$ and are not greater than $b$. The digraph $D$ is hamiltonian (is pancyclic, respectively) if it contains a hamiltonian cycle, i.e. a cycle of length $|V(D)|$ (contains a cycle of length $m$ for any $3\leq m \leq |V(D)|$).\\

Meyniel [12] proved the following theorem: if $D$ is a strong digraph on $n\geq 2$ vertices and $d(x)+d(y)\geq 2n-1$ for all pairs of non-adjacent vertices in $D$, then $D$ is hamiltonian (for short proofs of Meyniel's theorem see [5, 13]). 

Thomassen [15] (for $n=2k+1$) and Darbinyan [7] (for $n=2k$) proved: if $D$ is a digraph on $n\geq 5$ vertices with minimum degree at least $n-1$ and with minimum semi-degree at least $n/2-1$, then $D$ is hamiltonian (unless some extremal cases). \\

In each above mentioned theorems (as well as, in well known theorems Ghouila-Houri [10], Woodall [16], Manoussakis [11]) imposes a degree condition on all pairs of non-adjacent vertices (on all vertices). Bang-Jensen, Gutin, Li, Guo and Yeo [2, 3] obtained sufficient conditions for hamiltonisity of digraphs in which degree conditions requiring only for some pairs of non-adjacent vertices. Namely, they proved the following theorems (in all three theorems  $D$ is a strong digraph on $n\geq 2$ vertices).\\

\noindent\textbf{Theorem A} [2]. If $min \{d(x),d(y)\}\geq n-1$ and $d(x)+d(y)\geq 2n-1$ for every pair of non-adjacent vertices  $x$, $y$ with a common in-neighbour, then $D$ is hamiltonian.

\noindent\textbf{Theorem B} [2]. If $min \{d^+(x)+d^-(y),d^-(x)+d^+(y)\}\geq n$ for every pair of non-adjacent vertices $x$, $y$ with a common out-neighbour or a common in-neighbour, then $D$ is hamiltonian.

\noindent\textbf{Theorem C} [3]. If $min \{d^+(x)+d^-(y),d^-(x)+d^+(y)\}\geq n-1$ and $d(x)+d(y)\geq 2n-1$ for every pair of non-adjacent vertices $x$, $y$ with a common out-neighbour or a common in-neighbour, then $D$ is hamiltonian. \\
 Note that  Theorem C generalizes Theorem B.
In [9, 14, 6, 8] it was shown  that if a strong digraph $D$  satisfies the condition of the theorem of Ghouila-Houri [10] (Woodall [16], Meyniel [12], Thomassen and Darbinyan [15, 7]), then $D$ is pancyclic (unless some extremal cases, which are characterized). It is not difficult to check that the digraphs $ K^*_{n/2,n/2}$ and $ K^*_{n/2,n/2}-\{e\}$, where $n$ is even and $e$ is an arc of $ K^*_{n/2,n/2}$, satisfy the conditions of Theorem A (B, C) and has no cycle of odd length. Moreover, if in Theorems A, B, C  the digraph $D$ has no pair of non-adjacent vertices with a common in-neighbour and a common out-neighbour, then $D$ is a locally semicomplete digraph, and in [4], Bang-Jensen, Gutin and Volkmann characterize those strong locally semicomplete digraphs which are not pancyclic. For example, the following digraphs $D(5)$ and $D(6)$ with 5 and 6 vertices (respectively) are strong locally semicomplete, but has no cycle of length three, where
$$
V(D(5))=\{x_1,\ldots , x_4,y\} \,\, \hbox {and}\,\, A(D(5))=\{x_ix_{i+1}/1\leq i\leq 3\}\cup \{x_4x_1,x_2y,x_3y,yx_3,yx_4\};
$$
$$
V(D(6))=\{x_1,\ldots , x_5,y\} \,\, \hbox {and}\,\, A(D(6))=\{x_ix_{i+1}/1\leq i\leq 4\}\cup \{x_5x_1,x_1x_3,x_2x_4,x_2y,x_3y\} \cup
$$ 
$$ \{x_4y,yx_4,yx_5\}.
$$

It is natural to set the following problem:

\noindent\textbf{Problem}. Characterize those digraphs which satisfy the conditions of Theorem A (B, C), but are not pancylic.

To investigate that a given digraph $D$ is pancyclic, in [8, 13, 5, 7] it was proved the existence of cycles of length $|V(D)|-1$ or $|V(D)|-2$, and then using the constructions of these cycles it was proved that $D$ is pancyclic with some exceptions.

In this paper we prove two results which proide some support for the above problem: 

(i) if a strong digraph $D$ satisfies the condition of Theorem A and the minimum semi-degree of $D$ at least two; or 

(ii) if a strong digraph  $D$ is not directed cycle and satisfies the condition  of Theorem B, then either $D$ contains a cycle of length $n-1$ or $n$ is even and  $D$ is isomorphic to complete bipartite digraph or to complete bipartite digraph minus one arc.
  
Our proofs are based on the argument of [2, 3], which was turn based on the ideas used by  Bondy,  H\"{a}ggkvist and Thomassen  [5, 9, 14].  \\

\noindent\textbf{2. Preliminaries }\\

The following well-known simple lemmas is the basis of our results and other theorems on directed cycles and paths in digraphs. It we  will be used extensively in the proofs of our results.\\

\noindent\textbf{Lemma 1} [9]. Let $D$ be a digraph on $n\geq 3$ vertices containing a cycle $C_m$, $m\in [2,n-1] $. Let $x$ be a vertex not contained in this cycle. If $d(x,C_m)\geq m+1$, then  $D$ contains a cycle $C_k$ for all  $k\in [2,m+1]$.   \\

\noindent\textbf{Lemma 2} [5]. Let $D$ be a digraph on $n\geq 3$ vertices containing a path $P:=x_1x_2\ldots x_m$, $m\in [2,n-1]$ and let $x$ be a vertex not contained in this path. If one of the following conditions holds:

 (i) $d(x,P)\geq m+2$; 

 (ii) $d(x,P)\geq m+1$ and $xx_1\notin D$ or $x_mx_1\notin D$; 

 (iii) $d(x,P)\geq m$, $xx_1\notin D$ and $x_mx\notin D$;

\noindent\textbf{}then there is an  $i\in [1,m-1]$ such that $x_ix,xx_{i+1}\in D$, i.e. $D$ contains a path $x_1x_2\ldots x_ixx_{i+1}\ldots x_m$ of length $m$  (we say that  $x$ can be inserted into $P$ or the arc  $x_ix_{i+1}$ is a partner of $x$ on $P$). \\

\noindent\textbf{Lemma 3 }[2]. Let $P:=x_1x_2\ldots x_m$  be a path in $D$ and let $x$, $y$ be vertices of $V(D)-V(P)$ (possibly $x=y$). If there do not exist consecutive vertices $x_i, x_{i+1}$ on $P$ such that $x_ix$, $yx_{i+1}$ are arcs of $D$, then $d^-(x,P)+d^+(y,P)\leq m+\epsilon$, where $\epsilon =1$ if $x_mx\in D$ and 0, otherwise. \\

\noindent\textbf{3. Main results}\\

Let $C$ be a cycle in digraph $D$. For the cycle $C$, a $C$-bypass is an $(x,y)$-path $P$ of length at least two with both end-vertices $x$ and $y$ on $C$ and no other vertices on $C$. The length of the path $C[x,y]$ is the gap of $P$ with respect to $C$.

In the proof of Theorem 1, if $\{x,y\}$ is a pair of non-adjacent vertices with a common in-neighbour, then we say that $\{x,y\}$ is a good pair.\\

In the proofs of Theorems 1 and 2 we use (in the main) the notations which are used in the proofs of Theorems A and B (see [1], Theorems 5.6.1 and 5.6.5,  pages 248-250).\\

\noindent\textbf{Theorem 1}. Let $D$ be a strong digraph on $n$ vertices with minimum semi-degree at least two. Suppose that 
 $$
d(x)+d(y)\geq 2n-1 \quad \hbox {and} \quad min \{ d(x), d(y)\}\geq n-1 \eqno (*)
$$ 
for every pair of non-adjacent vertices $x$ and  $y$ with a common in-neighbour.
Then either $D$ contains  a cycle of length $n-1$ or $n$ is even and $D$ is  isomorphic to complete bipartite digraph $ K^*_{n/2,n/2}$ or $ K^*_{n/2,n/2}-\{e\}$, where $e$ is an arc of $ K^*_{n/2,n/2}$.

\noindent\textbf{Proof}. If $n\leq 4$, the theorem is easily verified. Suppose that $n\geq 5$ and $D$ contain no cycle of length $n-1$ and let $C:=x_1x_2\ldots x_mx_1$ be a longest non-hamiltonian cycle in $D$. Then  $3\leq m\leq n-2$ and let $R:=V(D)-V(C)$. Observe that if $y\notin V(C)$, then $y$ has no partner on $C$. We shall use this often without explicit reference.

We first prove the following claim: \\

\noindent\textbf{Claim 1}. Let $|R|\geq 3$ and $x_1yx_{\alpha +1}$ be a $C$-bypass of length two. If $A(y, C[x_2,x_{\alpha}])=\emptyset$, then $\alpha \geq 4$.

\noindent\textbf{Proof}. Suppose that $\alpha \leq 3$. Observe that $\{x_2, y\}$ is a good pair. Since $C$ is a longest non-hamiltonian cycle in $D$ and $|R|\geq 3$, it is not difficult to see that

$$
d^+(y,R)+d^-(x_i,R)\geq n-m-1  \quad \hbox {and } \quad d^-(y,R)+d^+(x_i,R)\leq n-m-1    \eqno (1)
$$
for every $i\in \{2,\alpha  \}$ and by Lemma 2(i),
$$
d(y,C)= d(y, C[x_{\alpha +1},x_1])\leq m-\alpha +2.     \eqno (2)
$$

 If $\alpha =2$, then  the vertex $x_2$ also cannot be inserted into $C[x_{\alpha +1},x_1]$ and hence by Lemma 2(i), $d(x_2)\leq m$. This together with (1) and (2) implies that $d(y)+d(x_2)\leq 2n-2$,  which is a contradiction since $\{y,x_2\}$ is a good pair.  So we can assume that $\alpha =3$. Now for $i\in \{2,3\}$ using (1) and (2) we obtain that
$$
d(y)+d(x_i)\leq 2n-m-3 +d(x_i,C).   \eqno (3)
$$
It is clear that if  $x_1x_i \in D$ ($i\in \{2,3\}$), then $\{y,x_i\}$ is a good pair and by (*) and (3) we obtain that $d(x_i,C)\geq m+2$. Therefore, by Lemma 2(i), $x_2$ has a partner on $C[x_4,x_1]$. This means  that $x_lx_2$, $x_2x_{l+1}\in D$ for some $l\in [4,m]$ and $x_3$ has no partner on $C[x_4,x_1]$ (otherwise, we obtain a  non-hamiltonian cycle longer than $C$). Now from above observation  it follows that $x_1x_3\notin D$. Observe that $x_3x_{l+1}\notin D$. If $a(x_{l+1},x_3)=0 $, then $\{x_3,x_{l+1}\}$ is a good pair and $d(x_3)\geq n-1$ by (*). It follows from (3) and $d(y)\geq n-1$ that $d(x_3,C[x_4,x_1])\geq m-1$ and hence by Lemma 2(ii), $x_1x_3\in D$, which is a contradiction. So we can assume that $x_{l+1}x_3 \in D$, $x_1\not= x_{l+1}$. Considering the pair $\{x_3, x_{l+2}\}$, we conclude analogously that $x_{l+2}x_3\in D$.  Continuing this process, we finally conclude  that $x_1x_3\in D$, contracting the conclusion above that this arc does not exist. Claim 1 is proved. \framebox  \\\\

In [2] (see [1], page 248), is proved that $D$ contains a $C$-bypass $P:=u_1u_2\ldots u_s$ ($s\geq 3$). W.l.o.g., let $u_1:=x_1$, $u_s:=x_{\gamma +1}$, $0< \gamma < m$. Suppose also that the gap $\gamma $ of $P$ is minimum among the gaps of all $C$-bypasses, i.e.
$$ A(\{u_2,\ldots , u_{s-1}\}, C[x_2,x_{\gamma}])=\emptyset.   \eqno (4) $$ 
Let $C':=C[x_2,x_{\gamma}]$ and $C'':=C[x_{\gamma +1},x_1]$. Note that if $\gamma \geq 2$, then $\{u_2,x_2\}$ is a good pair and by (*)
$$
min \{d(u_2),\, d(x_2)\}\geq n-1 \quad \hbox {and } \quad d(u_2)+ d(x_2)\geq 2n-1, \eqno (5)
$$
\\
We first show that $m\geq n-2$, i.e. $|R|\leq 2$. For this it suffices to consider the following four cases.

\noindent\textbf{Case 1}. $|R|\geq 3$ and $|R-P[u_2,u_{s-1}]|\geq 2$. The discussion of this case exactly is as same as the proof of Theorem 5.6.1 (see [1], page 249).\\ 

\noindent\textbf{Case 2}. $|R|\geq 3$ and $|R-P[u_2,u_{s-1}]|=1$. Then $\gamma \geq 3$ since $|P[u_2,u_{s-1}]|\geq 2$. Let $R-P[u_2,u_{s-1}]= \{y\}$. We can assume that $d^+(u_2, \{u_4,\ldots , u_{s-1}\})=0$ and $u_2x_{\gamma +1}\notin D$ (otherwise, we have Case 1). Therefore $d(u_2,R-\{y\})\leq |R|-1=n-m-1$, $d(u_2,R)\leq n-m+1$, and by Lemma 2(ii), $d(u_2,C'')\leq |C''|=m-\gamma +1$. This together with (4) and (5) implies that
 $$
n-1\leq d(u_2)= d(u_2,C'')+d(u_2,R-\{y\})+a(u_2,y)\leq n-\gamma +a(u_2,y).
$$
From this it is easy to see that $\gamma =3$,  $a(u_2,y)=2$, $d(u_2)=n-1$, $|R|=3$ since $\gamma \geq |R|$,  $d(u_2,C'')= n-5$ and $a(u_2,u_3)= 2$ ($s=4$). Then by (5), $d(x_2)\geq n$. Now since $d(x_2,C')\leq 2$ and $d(x_2,R)\leq 1$, it follows that $d(x_2,C'')\geq n-3\geq |C''|+2$. Therefore, by Lemma 2(i), $x_2$ has a partner on $C''$, i.e. $x_ix_2,x_2x_{i+1} \in D$, where $i\in [4,m]$. Thus the non-hamiltonian cycle $x_1u_2 u_{3}x_4\ldots x_ix_2x_{i+1}\ldots x_mx_1$ has length $n-2$, which is a contradiction.\\

\noindent\textbf{Case 3}. $|R|\geq 3$, $R=P[u_2,u_{s-1}]$ and $\gamma \not= 1$. Then $\gamma \geq |R|+1$, and $d(u_2,C'')\leq |C''|$ since $u_2x_{\gamma+1}\notin D$. We can assume that $d(u_2,R)\leq n-m$ (otherwise, we have Case 1 or 2). Therefore 
$$
n-1\leq d(u_2)=d(u_2,C'')+d(u_2,R)\leq |C''|+n-m\leq n-3
$$
since $|C''|= m-\gamma +1 $, which is a contradiction.

\noindent\textbf{Case 4}. $|R|\geq 3$, $R=P[u_2,u_{s-1}]$ and $\gamma =1$ (i.e. $x_{\gamma +1}=x_2$). We can assume that
if $2\leq i<j\leq s-1$, then $u_iu_j\in D$ if and only if $j=i+1$ (otherwise, we have one of the cases 1-3). Hence $d(u_2,R)\leq n-m$. Observe that 
$$
d^-(u_2,\{x_{m-1},x_m\})=0 \quad \hbox {and} \quad u_{s-1}x_3\notin D. \eqno (6)
$$

\noindent\textbf{Subcase 4.1}. $x_lu_2\in D$ for some $x_l\not= x_1$. Then by (6) there is a vertex $x_k$ with $2\leq k \leq m-2$ such that $x_ku_2\in D$ and $a(u_2,x_{k+1})=0$. Note that $x_k\rightarrow \{u_2,x_{k+1}\}$ and $\{u_2,x_{k+1}\}$ is a good pair. By (*), we have
$$
min \{d(u_2),\, d(x_{k+1})\}\geq n-1 \quad \hbox {and} \quad d(u_2)+ d(x_{k+1})\geq 2n-1 \eqno (7)
$$
Assume that $k$ is maximal with these properties. If $d^+(u_2, C[x_{k+2},x_1])=0$, then $A(u_2, C[x_{k+1},x_m])= \emptyset$ by the maximality of the $k$, and by Lemma 2(ii), $d(u_2,C)=d(u_2,C[x_1,x_{k}])\leq k$  since $u_2x_1\notin D$. Therefore $d(u_2)\leq k+n-m\leq n-2$, which is a contradiction. So we can assume that
 $d^+(u_2, C[x_{k+2},x_1])\not= 0$. Then there is an integer $\alpha \geq 1$, $k+\alpha \leq m$ such that $u_2x_{k+1+\alpha}\in D$  and $A(u_2,C[x_{k+1},x_{k+\alpha}])=\emptyset$. By Claim 1 we have $\alpha \geq 3$, and by Lemma 2(i), $d(u_2,C)\leq m-\alpha +1$. This along with $d(u_2,R)\leq n-m$ and (7) implies that
$$
n-1\leq d(u_2)= d(u_2,R)+d(u_2,C)\leq n-\alpha +1\leq n-2,
$$
which is a contradiction.

\noindent\textbf{Subcase 4.2}. $d^-(u_2, C[x_{2},x_m])=0$. Then $x_1\rightarrow \{ u_2,x_2 \}$, $a(x_2,u_2)=0$ and $\{ y,x_2 \}$ is a good pair. Therefore, by (*)
$$
 min \{d(u_2),\, d(x_{2})\}\geq n-1 \quad  \hbox {and} \quad d(u_2)+ d(x_{2}) \geq 2n-1 \eqno (8) 
$$ 
Note that $n-1\leq d(u_2)\leq n$. If $d(u_2)=n$, then it is not difficult to see that $d(u_2,C)=m$, $d(u_2,R)=n-m$, $u_2\rightarrow C[x_{3},x_1]$, $R-\{u_2\}\rightarrow y$ and $d^+(x_2, R)=d^-(x_2,R-\{u_{s-1}\})=0$. Therefore, since $x_2$ cannot be inserted into $C[x_3,x_1]$, we have $d(x_2)=d(x_2,R)+d(x_2,C)\leq m+1\leq n-2$, which  contradicts (8). So we can assume that $d(u_2)=n-1$ and by (8), $d(x_2)\geq n$.

If $uu_2\notin D$ for some vertex $u\in R-\{u_2\}$, then it is easy to see that $ R-\{u,u_2\}\rightarrow u_2$ and $u_2\rightarrow$ $C[x_3,x_1]$. Now we have 
$$d^-(x_2, R-\{u_{s-1}\})= d^+(x_2, R-\{u\})=0$$
and  hence, $d(x_2,R)\leq 2$. Therefore $n\leq d(x_2)=d(x_2,C)+d(x_2,R)\leq m+2\leq n-1$, which is a contradiction. Suppose that this is not the case, i.e.
 $R-\{u_2\}\rightarrow u_2$. Then $d^-(x_2, R-\{u_{s-1}\})=0$. If $yx_3\in D$, then it is easy to see that $d(x_2,C)\leq m$ and $d^+(x_2,R)=0 $. Therefore 
$d(x_2)\leq m+1 \leq n-2$, a contradiction. So we can assume that $u_2x_3\notin D$. Then $u_2\rightarrow C[x_4,x_1]$, $d(x_2,R)\leq 1$,  
 $A(x_3,\{u_2,u_{s-1}\})=\emptyset$  and  $d^-(u_{s-1},C[x_1,x_m])=0$
since $u_{s-1}u_2\in D$. Therefore $d^-(u_{s-1})=1$ since $d^-(u_{s-1},R)=0$, which contradicts that $d^-(u_{s-1}) \geq 2$.
 
Thus if $|R|\geq 3$, then in all possible cases we have obtained a contradiction. Therefore we have proved that $m=n-2$.  \framebox  \\\\

Let $R=\{y,z\}$. We first prove the following Claims 2-7.\\

\noindent\textbf{Claim 2}. If $x_iy, yx_{i+2}\in D$, $a(y,x_{i+1})=0$ and $d(x_{i+1})=n$, then $a(z,x_{i+1})=2$.
 
\noindent\textbf{Proof}. The proof of the claim immediately follows from the maximality of $C$ and Lemma 1. \framebox  \\\\

\noindent\textbf{Claim 3}. If $x_1y\in D$ and $a(y,x_2)=0$ (i.e., $\{x_2,y\}$ is a good pair). Then $a(y,x_3)\not=0$.

\noindent\textbf{Proof}. Suppose that the claim is not true, i.e.  $a(y,x_3)= 0$. By (*),

$$ min \{d(y),\, d(x_{2})\} \geq n-1 \quad  \hbox {and} \quad d(y)+ d(x_{2}) \geq 2n-1 \eqno (9) $$

Since $y$ cannot be inserted into $C[x_4,x_1]$, by Lemma 2(i) we have $d(y,C[x_4,x_1])\leq n-3$ (we can assume that $n\geq 6$). Therefore 
$$n-1\leq d(y)=a(y,z) +d(y,C[x_4,x_1])\leq n-1.$$ 
This implies that $a(y, z )=2$, $d(y,C[x_4,x_1])= n-3$ and $d(y)=n-1$. Therefore $d(x_2)\geq n$ (by (9)) and $yx_4\in D$ (by Lemma 2(ii)). Since $C$ is a longest non-hamiltonian cycle in $D$ and $a(y,z)=2$, $yx_4\in D$, it follows that $x_2z\notin D$, $zx_3\notin D$ and $d(x_2,C[x_4,x_1])\geq n-3$. 

Now we consider the following two possible cases.\\

\noindent\textbf{Case 1}. $a(x_2,z)=0$. Then $d(x_2,C[x_4,x_1])\geq n-2$ and by Lemma 2(i) $x_2$ has a partner on $C[x_4,x_1]$, $x_1z\notin D$ and $zx_4\notin D$. By Lemma 2(iii), $d(z,C[x_4,x_1])\leq n-5$. Therefore $d(z)\leq n-2$, since $d(z, \{y,x_3\})\leq 3$. This means that $z$ dose  not form a good pair with any vertex of $D$. Thus we have $d^-(z,C[x_3,x_1])=0$ since $x_1z\notin D$. Therefore $d^-(z)=1$, which is a contradiction.\\

\noindent\textbf{Case 2}. $zx_2\in D$. Then $x_2x_4\notin D$ and $x_my\notin D$ (otherwise, $C$ is not longest non-hamiltonian cycle in $D$). Since  
$d(x_2,C[x_4,x_1])\geq n-3$ and  $x_2x_4\notin D$, using Lemma 2(ii) we obtain, $x_2$ has a partner on 
$C[x_4,x_1]$, i.e. $x_ix_2,$ $x_2x_{i+1}\in D$ for some $i\in [4,m]$. Observe that $zx_4\notin D$ and $x_1z\notin D$, and by Lemma 2(iii), 
$
d(z,C[x_4,x_1])\leq n-5. 
$ 
Therefore $d(z)\leq n-1$ since $d(z,\{x_2,x_3\})\leq 2$.

If $a(x_3,z)=0$, then $d(z)\leq n-2$ and $z$  does not form a good pair with any vertex of $D$, which is not possible (since $yx_4\in D$, $zx_4\notin D$ and $x_1z\notin D$). Therefore $a(x_3,z)\not= 0$, i.e. $x_3z\in D$. It is easy to see that $yx_5\notin D$ and $x_my\notin D$. From this we obtain that $m\geq 5$. Now using Lemma 2(iii) we obtain, $d(y,C[x_5,x_m])\leq n-7$ and $yx_1$, $x_4y\in D$. Since $z$ has no partner on $C$ and $zx_2$, $x_3z\in D$, there is a vertex $x_l$ with $l\in [4,m+1]$ such that $x_{l-1}z\in D$ and $a(z,x_l)=0$ (i.e., $z$ forms a good pair with $x_l$). By (*) and $d(z)\leq n-1$ we have that $d(z)=n-1$. It follows that $zx_{l+1}\in D$ (by Lemma 2(i)). 
Let $l$ is minimal with these properties. Since $d(z)=n-1$, then $d(x_l)\geq n$ (by (*)), $d(x_l,C[x_{l+1},x_{l-1}])\leq n-2$ (by Lemma 2(i)) and $a(y,x_l)=2$ (Claim 2). Now we have, if $x_{l-2}z\in D$, then $x_{l-2}zyx_lx_{l+1}\ldots x_{l-2}$ is a cycle of length $n-1$. Therfore $x_{l-2}z\notin D$, i.e. $l=4$. This together with $x_2x_4\notin D$  and $d(x_4)\geq n$ implies that $x_4x_3\in D$ and $x_1yx_4x_3zx_5\ldots x_mx_1$ is a cycle of length $n-1$, a contradiction. Claim 3 is proved.  \framebox  \\\\ 

\noindent\textbf{Claim 4}. $d^-(y,\{x_i,x_{i+1}\})\leq 1$ for all $i\in [1,m]$.

\noindent\textbf{Proof}. Suppose, on the contrary, that (say) $\{x_m,x_1\}\rightarrow y$. W.l.o.g. we assume that  $a(y,x_2)=0$ (otherwise, $C\rightarrow y$, $d^+(y, C)=0$ and hence, $d^+(y)\leq 1$, a contradiction). This means that $\{y,x_2\}$ is a good pair. Therefore for $y$ and $x_2$ the condition (*) holds, i.e.
$$min \{d(y),\, d(x_{2})\}\geq n-1 \quad  \hbox {and} \quad d(y)+ d(x_{2}) \geq 2n-1.  $$
By Claim 3, we have that $a(y,x_3)\not=0 $. Since the minimum semi-degree of $D$ at least two, it is not difficult to see that $m\geq 4$. 
Consider the following three possible cases.\\

\noindent\textbf{Case 1}. $yx_3\in D$ and $yz\in D$. Then $d^+(z, \{x_2,x_3\})= 0$ and $d(x_2,C[x_3,x_1])\leq n-2$ (by Lemma 2(i)). Therefore $d(x_2)=n-1$ and $d(y)=n$ (by (*)), and $zy$, $x_2z\in D$, $x_1z\notin D$, $x_mx_2\notin D$.

First assume that $a(z,x_3)= 0$. Because of $x_2\rightarrow \{z,x_3\}$,  $\{z,x_3\}$ is a good pair and  by(*), $d(z),\, d(x_3)\geq n-1$. Since $x_1z\notin D$, using Lemma 2(ii), we obtain that $zx_4\in D$ and $d(z)=n-1$ (otherwise, $d(z)=d(z,\{y,x_2\})+d(z,C[x_4,x_1])\leq n-2$, a contradiction). Hence $d(x_3)\geq n$ by (*), and $x_3y\in D$ by Claim 2. It is easy to see that $x_1x_3\notin D$, and by Lemma 2(ii), $x_3x_2\in D$ since $d(x_3)\geq n$ and $x_3$ cannot be inserted into $C[x_4,x_1]$. Therefore $x_myx_3x_2zx_4\ldots x_m$ is a cycle of length $n-1$, a contradiction.

Second assume that $x_3z\in D$. Then since $\{x_2,x_3\}\rightarrow z$ and $x_1z\notin D$ there is an integer $k\in[3,m]$ such that $\{x_{k-1},x_k\}\rightarrow z$ and $a(z,x_{k+1})=0$, i.e. $\{x_{k+1},z\}$ is a good pair. Since $zx_2\notin D$ and $x_1z\notin D$ using Lemma 2(ii) we obtain, 
$$
d(z)=d(z,C[x_2,x_k])+d(z,C[x_{k+2},x_1])+a(z,y)\leq n-1.
$$
This together with (*) and $d^+(z)\geq 2$ implies that $d(x_{k+1})\geq n$ and $zx_{k+2}\in D$. Therefore $a(y,x_{k+1})=2$ (Claim 2) and $x_1x_2\ldots x_{k-1}zyx_{k+1}\ldots x_mx_1$ is a cycle of length $n-1$, a contradiction.\\

\noindent\textbf{Case 2}. $yx_3\in D$ and $yz\notin D$.
From $d(y)\geq n-1$ and $d(y,C[x_3,x_1])\leq n-2$ it follows that $zy\in D$ and  $d(y)= n-1$. Then $d(x_2)\geq n$ (by (*)) and $a(x_2,z)=2$ (Claim 2). Hence $yx_4\notin D$ and by Lemma 2(ii), $d(y,C[x_4,x_1])\leq n-4$.  This along  with $yz\notin D$ and $d(y)=n-1$ implies that $x_3y\in D$. Observe that $x_mx_2\notin D$ and  by Lemma 2(ii), $d(x_2,C[x_3,x_m])\leq n-4$ since $x_2$ cannot be inserted into $C[x_3,x_m]$. Therefore $x_2x_1\in D$.

\noindent\textbf{Subcase 2.1}. $a(z,x_3)=0$. Then $\{z,x_3\}$ is a good pair. This together with (*), $yz\notin D$ and Lemma 2(i) implies that $d(z)=n-1$, $zx_4\in D$ and $d(x_3)\geq n$. It is clear that $x_3x_2\notin D$ (otherwise, $C_{n-1}:=x_myx_3x_2zx_4\ldots x_m$). Now again using Lemma 2(i), we obtain that $x_1x_3\in D$ since $d(x_3,C[x_4,x_1])\geq n-3$. Observe that $x_mz\notin D$ (otherwise, $C_{n-1}:=x_mzx_2x_1x_3\ldots x_m$) and
$$
n-1=d(z)= d(z,C[x_4,x_m])+d(z,\{y,x_1,x_2,x_3\})\leq n-1.
$$
From this it follows that $zx_1\in D$. Let $m\geq 5$. Then $x_3x_5\notin D$ (otherwise, $C_{n-1}:=x_1x_2zyx_3x_5\ldots x_mx_1$) and by Lemma 2(ii), $x_4x_3\in D$ since $d(x_3)=n$. Thus we have a cycle $C_{n-2}:=x_1x_2zx_4\ldots x_mx_1$  which does not contain the vertices $y,x_3$ and $\{x_1,x_2\}\rightarrow x_3$, $a(z,x_3)=0$ and $x_3y$, $x_3x_4\in D$. Therefore for this cycle $C_{n-2}$ Case 1 holds.

Let now $m=4$, i.e. $n=6$. Then because of $d^-(z,\{y,x_1,x_3,x_m\})=0$ we have $d^-(z)=1$, which is a contradiction. \\

\noindent\textbf{Subcase 2.2}. $a(z,x_3)\not= 0$. Then $x_3z\in D$ since $a(x_2,z)=2$. Note that $x_{m-1}z\notin D$ (otherwise, $C_{n-1}:=x_{m-1}zx_2x_1yx_3$ $\ldots x_{m-1}$). Then $m\geq 5$ and it is easy to see that there is an integer $k\in [3,m-2]$ so that $\{x_{k-1},x_k\}\rightarrow z$, $a(z,x_{k+1})=0$. Then $\{z,x_{k+1}\}$ is a good pair since $x_k\rightarrow \{z,x_{k+1}\}$. From this by (*), $d(z)\geq n-1$, and by Lemma 2 (i), $zx_{k+2}\in D$. Therefore for the vertex $z$ we have considered Case 1. \\

\noindent\textbf{Case 3}. $yx_3\notin D$. Then $x_3y\in D$ by Claim 3. Using Lemma 2(ii) and (*), it is not difficult to see that 
$$
d(y,C)=n-3,\,\, a(y,z)=2,\,\,d(y)=n-1,\,\, d(x_2)\geq n \, \,  \hbox {and}\,\, d^+(z, \{x_2,x_3,x_5\})=0. $$

If $x_4y\in D$, then $m\geq 6$ and from  $d^+(y)\geq 2$ and $d(y)=n-1$ it follows that  $\{x_{k-1},x_k\}\rightarrow y \rightarrow \{x_{k+2}\}$ and $a(y,x_{k+1})=0$ for some $k\in [4,m-2]$, i.e. we have the considered Case 1.

So we can assume that  $x_4y\notin D$. Then $m\geq 5$,  $a(y,x_4)=0$ and $\{y,x_4\}$ is a good pair since $x_3\rightarrow \{x_4,y\}$. Therefore by (*), $d(x_2)$ and $d(x_4)\geq n$ since $d(y)=n-1$. From $d(y, \{x_2,x_3,x_4, z)\})=3$ it follows that  $d(y,C[x_5,x_1])=n-4$. Hence $yx_5\in D$ (by Lemma 2(ii)) and $a(x_4,z)=2$ (Claim 2). Now it is easy to see that  $x_2x_4\notin D$.  Again using Lemma 2(ii), we obtain that $d(x_2,C[x_4,x_1])\leq n-4$ and $d(x_2)\leq n-1$ since $zx_2\notin D$, contradicting the conclusion above that $d(x_2)\geq n$. 
 Claim 4 is proved. \framebox  \\\\

\noindent\textbf{Claim 5}. If $x_iy\in D$ and $a(y,x_{i+1})=0$, then $d^+(y,\{x_{i+2},x_{i+3}\})\leq 1$ for all $i\in[1,m]$.

\noindent\textbf{Proof}. Suppose that the claim is not true. W.l.o.g we can assume that $x_1y\in D$, $a(y,x_2)=0$ and $y\rightarrow \{x_3,x_{4}\}$. Note that $\{y,x_2\}$ is a good pair.

\noindent\textbf{Case 1}. $zy\notin D$. Then from $d(y,C)\leq n-2$ and (*) it follows that $yz\in D$, $d(y)=n-1$ and $d(x_2)\geq n$. Observe that $a(x_2,z)=2$ (Claim 2) and $a(z,x_3)=0$ (Claim 4). Therefore $\{z,x_3\}$ also is a good pair and for its (*) holds. Using Lemma 2(i) and (*), it is not difficult to see that $d(z)=n-1$, $d(x_3)\geq n$ and $zx_4\notin D$ since $zy\notin D$. Now by Claim 2, $a(x_3,y)=2$, which is a contradiction.

\noindent\textbf{Case 2}. $zy\in D$. Then $d^-(z,\{x_1,x_2\})=0$, $d(x_2,C)=n-2$, $zx_2\in D$, $d(y)=n$, $yz\in D$ and $zx_3\notin D$ (i.e., either $x_3z\in D$ or $a(z,x_3)=0$) by Lemma 2(i) and (*).

\noindent\textbf{Subcase 2.1}. $a(z,x_3)\not=0$. Then  it is easy to see that $x_3z\in D$, $m\geq 4$ and $yx_5\notin D$, $d(y,C[x_5,x_1])=n-5$, $x_4y\in D$ since $d(y)=n$. By Claim 4,  $x_4z\notin D$. Therefore $a(z,x_4)=0$ and $\{z,x_4\}$ is a good pair and hence, $d(z)$, $d(x_4)\geq n-1$. Using Lemma 2(ii) and (*) we obtain that $zx_5\in D$, $d(z)=n-1$, $d(x_4)=n$  since $x_2z\notin D$. If $x_2x_4\in D$, then $C_{n-1}:=x_2x_4yzx_5\ldots x_mx_1x_2$, and if $x_2x_4\notin D$, we obtain that $x_4x_3\in D$ (since $x_4$ cannot be inserted into $C[x_5,x_3]$ and $d(x_4)=n$) and $C_{n-1}:=x_1yx_4x_3zx_5\ldots x_mx_1$, which is a contradiction.

\noindent\textbf{Subcase 2.2}. $a(z,x_3)=0$. Then $\{z,x_3\}$ is a good pair since $y\rightarrow \{z,x_3\}$, and $m\geq 4$ since $d^-(z)\geq 2$. Therfore by (*), $x_2z\notin D$ and Lemma 2(ii), we obtain that $zx_4\in D$, $d(z)=n-1$ and $d(x_3)\geq n$. Since the vertex $x_3$ has no partner on the cycle $C_{n-2}:=x_1yzx_4\ldots x_mx_1$ and $d(x_3)\geq n$, using Lemma 2(i) we obtain that $x_1x_3\in D$. Now for this cycle $C_{n-2}$ we have $\{x_1,y\}\rightarrow x_3$, which contradicts Claim 4. Claim 5 is proved. \framebox \\\\

\noindent\textbf{Claim 6}. If $a(y,z)=1$, then $n$ is even and $D\equiv K^*_{n/2,n/2}-\{e\}$, where $e$ is an arc of $K^*_{n/2,n/2}$.

\noindent\textbf{Proof}. Since $D$ contain no cycle of length $n-1$, using Lemma 1 we obtain  that $d(y)$, $d(z)\leq n-1$. W.l.o.g. assume that $yz$,  $x_my\in D$, $a(y,x_1)=0$ (if $C\rightarrow y$, then $D$ contains a cycle of length $n-1$ since  $d^+(z, C)\not=0$). Then $\{x_1,y\}$ is a good pair and hence, $d(y)=n-1$, $yx_2\in D$ and $d(x_1)\geq n$ by Lemma 2(i) and (*). We have $a(x_1,z)=2$ (Claim 2), and $d(x_1,C[x_2,x_m])=n-2$ by  Lemma 2(i). Therefore $d^-(z,\{x_m,x_2\})=0$ (Claim 4) and $a(z,x_2)=0$, i.e. $\{z,x_2\}$ is a good pair, $d(x_2)\geq n$ and  by Lemma 2(ii), $d(z,C[x_3,x_m])=n-4$, $zx_3\in D$ (since $x_mz\notin D$). From this $a(x_2,y)=2$ (Claim 2), $zx_4\notin D$, $a(y,x_3)=0$ (Claim 4), $d(x_3)\geq n$ (by (*)) and $d(y,C[x_4,x_m])= n-4$. Therefore $yx_4\in D$, $x_3z\in D$ (Claim 2), $a(z,x_4)=0$ (Claim 4), $d(x_4)\geq n$ (by (*)) and $zx_5\in D$. Continuing this process, we finally conclude that $n$ is even, $d(x_i)=n:=2k$, 
$$ y\rightarrow \{x_2,x_4,\ldots , x_{2k-2}\} \rightarrow y, \quad z \rightarrow \{x_1,x_3,\ldots , x_{2k-3}\}\rightarrow z $$ and
$$
A(y,\{x_1,x_3,\ldots , x_{2k-3}\})=A(z,\{x_2,x_4,\ldots , x_{2k-2}\})=\emptyset.
$$

 Now we prove that 
$$
A(\langle\{x_1,x_3,\ldots , x_{2k-3}\}\rangle)=A(\langle\{x_2,x_4,\ldots , x_{2k-2}\}\rangle)=\emptyset.
$$
 
Suppose this is not the case.
Let $x_ix_j\in D$, where $i,j\in \{1,3,\ldots , 2k-3 \}$. Then 
$$ 
a(x_i,z)=a(x_j,z)=a(x_{i-1},y)=a(x_{i+1},y)=a(x_{j-1},y)=a(x_{j+1},y)=2.
$$
If $|C[x_i,x_j]|=3$, then $C_{n-1}:=x_ix_j\ldots x_{i-1}yzx_i$; if $|C[x_i,x_j]|\geq 5$, then $C_{n-1}:=x_ix_j\ldots x_{i-1}yx_{i+1}\ldots$ $x_{j-2}zx_i$.
Let now  $x_ix_j\in D$, where $i,j\in \{2,4,\ldots , 2k-2 \}$. Then 
$$ 
a(x_i,y)=a(x_j,y)=a(x_{i-1},z)=a(x_{i+1},z)=a(x_{j-1},z)=a(x_{j+1},z)=2.
$$
If $|C[x_i,x_j]|=3$, then $C_{n-1}:=x_ix_j\ldots x_{i-2}yzx_{i-1}x_i$; if $|C[x_i,x_j]|\geq 5$, then $m\geq 6$ and $C_{n-1}:=x_ix_j\ldots x_{i-1}zx_{i+1}$ $\ldots  x_{j-2}yx_i$. In all possible cases we have that $D$ contains a cycle of length $n-1$, which is a contradiction. Therefore 
$$
A(\langle\{y,x_1,x_3,\ldots , x_{2k-3}\}\rangle)=A(\langle\{z,x_2,x_4,\ldots , x_{2k-2}\}\rangle)=\emptyset,
$$
i.e. $D\equiv K^*_{n/2,n/2}-\{e\}$.
Claim 6 is proved. \framebox \\\\

\noindent\textbf{Claim 7}. If $a(y,z)=2$ and $d(y)=n$, then $n$ is even and either $D\equiv K^*_{n/2,n/2}$ or $D\equiv K^*_{n/2,n/2}-\{e\}$, where $e$ is an arc of $K^*_{n/2,n/2}$. 

 \noindent\textbf{Proof}. For definite let $x_my\in D$ and $a(y,x_1)=0$. Then $\{y,x_1\}$ is a good pair. By Claim 3, $a(y,x_2)\not=0$, and  since  $d(y,C)= n-2$, using Lemma 2(i) we obtain that $yx_2\in D$. It is not difficult to see that $m\geq 4$, and by Claim 4, $x_{m-1}y\notin D$. From $x_my, yx_2\in D$ and Claim 5 it follows that $yx_3 \notin D$. Now by Lemma 2(ii) we have $d(y,C[x_3,x_m])\leq n-4$. Therefore $x_2y\in D$ and by Claim 4, $a(x_3,y)=0$ and hence, $d(y,C[x_4,x_m])= n-4$. Again using  Lemma 2(ii) we obtain that $yx_4\in D$ and hence by Claim 5, $yx_5\notin D$. Therefore, by Lemma 2(ii), $d(y,C[x_5,x_m])\leq n-6$ and hence, $x_4y\in D$ and  $a(y,x_5)=0$. Continuing this process we finally conclude that $n:=2k$ and 
 $$ 
a(y,x_2)=a(y,x_4)=\ldots =a(y,x_{2k-2})=2, \quad a(y,x_1)=a(y,x_3)=\ldots =a(y,x_{2k-3})=0.
$$
 Observe that $d(x_{2i-1})\geq n-1$ by (*) for every $i\in [1,k-1]$. Consider the cycle $C_{n-2}:=x_{2i}yx_{2i+2}\ldots x_{2i}$ of length $n-2$, $i\in [1,k-1]$. Note that the vertices  $x_{2i+1}$ and $z$ are not on this cycle. We can assume that $a(z,x_{2i+1})=2$ (otherwise, by Claim 6, $D\equiv K^*_{n/2,n/2}-\{e\}$). Analogously to the proof of Claim 6, we get that
 $$
A(\langle\{y,x_1,x_3,\ldots , x_{2k-3}\}\rangle)=A(\langle\{z,x_2,x_4,\ldots , x_{2k-2}\}\rangle)=\emptyset.
$$
Now using the condition (*) and the fact that for every pair of distinct $i,j\in \{1,3,\ldots ,2k-3\}$ ($i,j\in \{2,4,\ldots ,2k-2\}$), $\{x_i,x_j\}$ is a good pair, we conclude that either $D\equiv K^*_{n/2,n/2}$ or $D\equiv K^*_{n/2,n/2}-\{e\}$. Claim 7 is proved. \framebox \\\\

   Let us now complete the proof of the theorem. By Claims 6 and 7 we can assume that for any cycle of length $n-2$ in $D$ if the vertices $u$ and $v$ are not on this cycle then $max \{d(u),\, d(v)\}\leq n-1$ and  $a(u,v)=2$. 

W.l.o.g. assume that $x_my\in D$ and $a(y,x_1)=0$, i.e. $\{y,x_1\}$ is a good pair. Then $d(y)=n-1$ and $d(x_1)=n$ by the  our assumption and (*).  Then $a(y,x_2)\not=0$ by Claim 3. Let  $yx_2\in D$, then  $C_{n-2}:=x_myx_2\ldots x_m$ and $d(x_1)=n$, which contradicts to our assumption. Let now  $yx_2\notin D$. Then $x_2y\in D$ and since $d^+(y)\geq 2$ and $d(y,C[x_2,x_m])=n-3$, it is not difficult to see that for some $j\in [2,m-2]$, $x_jy,yx_{j+2}\in D$ and $a(y,x_{j+1})=0$. A similar argument applies for this case, we again obtain a contradiction. The proof of Theorem 1 is complete. \framebox  \\\\

The following example shows that the sharpness  the  minimum semi-degree condition in Theorem 1 would be best possible in the sense that for all $n=k+2 \geq 6$ there is a strong digraph $D$ on $n$ vertices which has minimum semi-degree one and satisfies the condition (*) of Theorem 1, but contain no cycle of length $n-1$. To see this,
let $D$ be a digraph with vertex set $V(D)=\{y,z, x_1,x_2,\ldots , x_k\}$; and let (for the convenience of the reader) $N^-(y)=\{z,x_1,x_3,x_4,\ldots , x_k\}$ and $N^+(y)=\{z\}$;  $N^-(z)=\{y,x_1,x_2,x_4,x_5,\ldots , x_k\}$ and $N^+(z)=\{y, x_4\}$; $N^-(x_1)=\{x_k,x_2,x_3\}$ and $N^+(x_1)=\{y,z,x_2, x_4\}\cup \{x_5,x_6,\ldots , x_{k-1}\}$; $N^-(x_2)=\{x_1,x_3\}$ and $N^+(x_2)=\{z,x_1, x_3,x_4,\ldots , x_{k}\}$; $N^-(x_3)=\{x_2\}$ and $N^+(x_3)$ $=\{y,x_1, x_2,x_4,x_5,\ldots , x_{k}\}$; $N^-(x_4)=\{z,x_1,x_2,x_3\}\cup \{x_6, x_7,\ldots , x_k\}$ and  $N^+(x_4)=\{y,z,x_1\}$ if $k=4$ and $N^+(x_4)=\{y,z,x_5\}$ if $k\geq 5$;  if $5\leq i \leq k-1$, then $N^-(x_i)=\{x_1,x_2,x_3,x_{i-1}\}\cup \{x_{i+2},x_{i+3},\ldots , x_k\}$ and $N^+(x_i)=\{y,z,x_{i+1}\}\cup \{x_4,x_5,\ldots , x_{i-2}\}$; finally if $k\geq 5$, then let $N^-(x_k)=\{x_2,x_3,x_{k-1}\}$ and $N^+(x_k)=\{y,z,x_1\}\cup \{x_4,x_5,\ldots , x_{k-2}\}$, where $\{x_{i},x_{i+1},\ldots , x_j\}=\emptyset$ if $j\leq i-1$. 

Note that $x_1x_2\ldots x_kx_1$ is a cycle of length $k=n-2$, $\langle \{x_1,x_2,\ldots , x_k\}\rangle$ is a semicomplete digraph, the pairs of non-adjacent distinct vertices with a common in-neighbour are only $\{y,x_2\}$ and  $\{z,x_3\}$, $d(y)=d(x_3)=n-1$, $d(z)=d(x_2)=n$ and $d^+(y)=d^-(x_3)=1$. It is not difficult to check that $D$ is strong, satisfies the condition (*) of Theorem 1 and contain no cycle of length $n-1$. 

Moreover the following example from [14] (also [1], p. 300) also shows that in Theorem 1 the minimum semi-degree condition ($\geq 2$) cannot be replaced by one. For some $m\leq n$ let $D_{n,m}$ be the digraph with vertices $V(D_{n,m})=\{x_1,x_2,\ldots , x_n\}$ and arcs $A(D_{n,m})=\{x_ix_j/ i<j $ or $i=j+1\}$ $\setminus \{x_ix_{i+m-1}/1\leq i\leq n-m+1\}$. $D_{n,m}$  is strong, has no cycle of length $m$ and if $m=n-1$, then the pairs $\{x_1,x_{n-1}\}$ and  $\{x_2,x_{n}\}$ are only the non-adjacent pairs with a common in-neighbour. It is easy to check that $d^-(x_1)=d^+(x_n)=1$, $d(x_1)=d(x_n)=n-1$ and $d(x_{n-1})=d(x_2)=n$.\\

\noindent\textbf{Theorem 2}. Let $D$ be a strong digraph on $n\geq 4$ vertices, which is not directed cycle of length $n$. Suppose that 
 $$min \{d^+(x)+d^-(y),\, d^-(x)+d^+(y)\}\geq n \eqno (**)$$ 
for every pair of non-adjacent vertices $\{x, y\}$ with a common out-neighbour or a common in-neighbour. Then either $D$ contains  a cycle of length $n-1$ or $n$ is even and $D$ isomorphic to complete bipartite digraph $ K^*_{n/2,n/2}$. 

\noindent\textbf{Proof}. Suppose that $D$ has no cycle of length $n-1$ and $C:=x_1x_2\ldots x_mx_1$ is a longest non-hamiltonian cycle in D. Let $R:=V(D)-V(C)$. Then $3\leq m \leq n-2$, i.e. $|R|\geq 2$. In [2] (see [1] page 250), was shown that $D$ has a $C$-bypass  with three vertices. W.l.o.g. assume that $B:=x_1yx_{j+1}$  is a $C$-bypass and the gap $j$ of $B$ with respect to $C$ is minimum among the gaps of all $C$-bypasses with three vertices. Clearly, $j\geq 2$ and
$$
A(y,\{x_2,x_3,\ldots  , x_j\})=\emptyset.   \eqno (10)
$$
Observe that $\{y,x_2\}$ ($\{y,x_j\}$, respectively) is a pair of non-adjacent vertices with a common in-neighbour $x_1$(with a common out-neighbour $x_{j+1}$, respectively). Therefore for these pairs the condition (**) of the theorem holds.

Let $C'':=C[x_{j+1},x_1]$ and $C':=C[x_{2},x_j]$. By Lemmas 2, 3 and the maximality of $C$ we have 
$$
d(y,C'')\leq |C''|+1; \quad d^-(x_2,C'')+d^+(x_j,C'')\leq |C''|+1. \eqno (11)
$$

\noindent\textbf{Case 1}. $|R|\geq 3$, i.e. $m\leq n-3$. Then
$$
d^+(y,R)+d^-(x_2,R)\leq n-m-1\quad \hbox{and} \quad d^+(x_j,R)+d^-(y,R)\leq n-m-1,
$$
otherwise, $D$ contains a long non-hamiltonian cycle than $C$.  This along with (10), (11) and (**) gives 
$$
2n\leq d^-(y)+d^+(x_j)+ d^+(y)+d^-(x_2)= d^-(y,R)+d^+(x_j,R)+ d^+(y,R)+d^-(x_2,R)+ d(y,C'')+$$
 $$ d^+(x_j,C'') 
+ d^-(x_2,C'')+d^-(x_2,C')+d^+(x_j,C')\leq 2(n-m-1)+2|C''|+2+2|C'|-2\leq 2n-2 $$
since $|C'|+|C''|=m$, which is a contradiction.\\

\noindent\textbf{Case 2}. $|R|=2$, i.e. $m=n-2$. Let $R=\{y,z\}$.

\noindent\textbf{Subcase 2.1}. $j=2$, i.e. $x_1y,yx_3\in D$ and $a(x_2,y)=0$. Then $|C''|=n-3$. By the maximality of the cycle $C$ and Lemma 1 we have $d(y,C), d(z,C)\leq n-2$. From the condition (**) of the theorem it follows that $d(y)\geq n$ or $d(x_2)\geq n$. W.l.o.g. we assume that $d(y)\geq n$. Then 
$$n\leq d(y)=a(y,z)+d(y,C)\leq 2+d(y,C).$$ 
Since $d(y,C)\leq n-2$, it follows that $a(y,z)=2$, $d(y,C)=n-2$ and $d(y)=n$. Similarly, we obtain that have $d(x_2)=n$, $d(x_2,C)=n-2$ and $a(x_2,z)=2$ (by (**) and (11)). It is easy to see that $n\geq 6$. Observe that $yx_4\notin D$ and $x_my\notin D$. Now using Lemma 2(iii) we obtain, $d(y,C[x_4,x_m])=n-6$ and $a(x_1,y)=a(y,x_3)=2$. If $n=6$, then it is easy to cheek that $D\equiv K^*_{3,3}$. Assume that $n\geq 7$. If $\{x_{k-1},x_k\}\rightarrow y$ and $a(y,x_{k+1})=0$ for some $k\in [4,m-1]$, then by Lemma 2(ii), $yx_{k+2}\in D$. Then, since $\{y,x_{k+1}\}$ is a pair of non-adjacent vertices with a common in-neighbour $x_k$, $d(y)=n$ and the vertex $x_{k+1}$ has no partner on $C[x_{k+2},x_k]$ it follows that $d(x_{k+1})=n$ and $a(x_{k+1},z)=2$. Therefore $C_{n-1}:=x_{k+1}C[x_{k+2},x_{k-1}]yzx_{k+1}$, a contradiction. So we can assume that $d^-(y,\{x_i,x_{i+1}\})\leq 1$ for all $i\in [1,m]$. This together with $x_3y\in D$ and $yx_4\notin D$ implies that $a(x_4,y)=0$. Analogously above, we obtain that $d^+(y, \{x_i,x_{i+1}\})\leq 1$ for all $i\in [1,m]$. Now it is not difficult to see that $n$ is even ($n:=2k+2$), $a(y,x_i)=2$ for all $i\in \{1,3,\ldots , 2k-1\}$ and $A(y,\{x_2,x_4,\ldots ,x_{2k}\})=\emptyset$. Then $\{y,x_{2j}\}$ is a pair of non-adjacent vertices with a common in-neighbour $x_{2j-1}$ for all $j\in [1,k]$. Therefore $d(x_{2j})=n$ since $d(y)=n $, and $a(z,x_{2j})=2$, $A(z,\{x_1,x_3,\ldots ,x_{2k-1}\})=\emptyset$ since $x_{2j}$ cannot be inserted into $C[x_{2j+1},x_{2j-1}]$. We finally conclude that either $D$ contains a cycle of length $n-1$ or $D\equiv K^*_{n/2,n/2}$ with partite sets $\{z,x_1,x_3,\ldots , x_{2k-1}\}$ and $\{y,x_2,x_4,\ldots , x_{2k}\}$.

\noindent\textbf{Subcase 2.2}. $j\geq 3$. From (11) and (**) it follows that
$$
min\{d(x_2), d(x_j)\}\geq 2n-d(y)= 2n-d(y,C'')-a(y,z)\geq  n+j-a(y,z). $$
Therefore  
$$
d(x_2,C'')\geq n+j-2(|C'|-1)-d(z,\{y,x_2\})\geq  n-j+4-d(z,\{y,x_2\}) 
$$
and similarly
$$
 d(x_j,C'')\geq  n-j+4-d(z,\{y,x_j\}). 
$$

From this it follows that if $j=3$, then $d(z,\{y,x_2\})$,  $d(z,\{y,x_2\})\leq 3$ and $d(x_2,C'')$,  $d(x_j,C'')\geq n-j+1=|C''|+2$. So, by Lemma 2(i) we have that $x_2$ and $x_3$ has a partner on $C''$ and therefore, $D$ contains a cycle of length $n-1$, a contradiction. Now we can assume that $j\geq 4$. Note that $d(y,C)\leq n-j$ by Lemma 2(i).

First assume that $a(y,z)\leq 1$. Then $d(y)\leq n-j+1$ and by (**)
$$
2n\leq d^+(y)+d^-(x_2)+d^-(y)+d^+(x_j)\leq n-j+1+d^-(x_2)+d^+(x_j) 
$$
and
$$
d^-(x_2)+d^+(x_j)\geq n+j-1.\eqno (12)
$$
This together with (11) implies that 
$$
n+j-1\leq d^-(x_2)+d^+(x_j)\leq n-j+2(j-2)+2=n+j-2,
$$
which is a contradiction.

Second assume that $a(y,z)=2$. Then $d(y)\leq n-j+2$ and similarly (12) we obtain that $d^-(x_2)+d^+(x_j)\geq n+j-2$. On the other hand using (11) it is easy to see that $d^-(x_2)+d^+(x_j)\leq n+j-2$. Therefore $d^-(x_2)+d^+(x_j)= n+j-2$, $zx_2$, $x_jz\in D$ and $d(y, C'')=n-j\geq 3$. From this it is not difficult to see that $x_{m}y\notin D$ and $yx_{j+2}\notin D$. Therefore $n-j\geq 4$, and by Lemma 2(iii), $d(y,C[x_{j+2},x_{m}])\leq n-j-4$. Hence $yx_1, x_{j+1}y\in D$ and $d(y,C[x_{j+2},x_{m}])= n-j-4$. Now using Lemma 2 we obtain that $x_{i-1}y,yx_{i+1}\in D$  and $a(x_i,y)=0$ for some $i\in [j+2,m]$, i.e. we have the considered Subcase 2.1. The theorem is proved. \framebox  \\\\

We believe Theorem 2 can be generalized to the following

\noindent\textbf {Conjecture}. Let $D$ be a strong digraph on $n\geq 4$ vertices. Suppose that $min \{d^+(x)+d^-(y),d^-(x)+d^+(y)\}\geq n-1$ and $d(x)+d(y)\geq 2n-1$ for every pair of non-adjacent vertices $x$, $y$ with a common out-neighbour or a common in-neighbour (i.e., satisfies the conditions of Theorem C). Then $D$ contains a cycle of length $n-1$ maybe except some digraphs which has a "simple" characterization.\\

\noindent\textbf {References}\\

[1] J. Bang-Jensen, G. Gutin, Digraphs: Theory, Algorithms and Applications, Springer, 2000.

[2] J. Bang-Jensen, G. Gutin, H. Li, Sufficient conditions for a digraph to be hamiltonian, J. Graph Theory 22 (2) (1996) 181-187.

[3] J. Bang-Jensen, Y. Guo, A.Yeo, A new sufficient condition for a digraph to be hamiltonian, Discrete Applied Math., 95 (1999) 77-87.
 
[4] J. Bang-Jensen, Y. Guo, L. Volkmann, A classification of locally semicomplete digraphs. 15th British Combinatorial Conference (Stirling, 1995). Discrete Math.. 167/168 (1997) 101-114.

[5] J.A. Bondy, C. Thomassen, A short proof of Meyniel's theorem, Discrete Math. 19 (1977) 195-197.

[6] S.Kh. Darbinyan,  Pancyclic and panconnected digraphs, Ph. D. Thesis, Institute  Mathematici Akad. Nauk BSSR, Minsk, 1981 (see also, Pancyclicity of digraphs with the Meyniel condition, Studia Sci. Math. Hungar., 20 (1-4) (1985) 95-117, in Russian).

[7] S.Kh. Darbinyan, A sufficient condition for the Hamiltonian property of digraphs with  large semidegrees, Akad. Nauk Armyan. SSR Dokl. 82 (1) (1986) 6-8 (see also, arXiv: 1111.1843v1 [math.CO] 8 Nov 2011).

[8] S.Kh. Darbinyan, On the  pancyclicity of digraphs with large semidegrees,  Akad. Nauk Armyan. SSR Dokl. 83 (3) (1986) 99-101 (see also, arXiv: 1111.1841v1 [math.CO] 8 Nov 2011).

[9] R. H\"{a}ggkvist, C. Thomassen, On pancyclic digraphs, J. Combin. Theory Ser. B 20 (1976) 20-40.

[10] A. Ghouila-Houri, Une condition suffisante d'existence d'un circuit hamiltonien, C. R. Acad. Sci. Paris Ser. A-B 251 (1960) 495-497.

[11] Y. Manoussakis, Directed Hamiltonian graphs, J. Graph Theory 16 (1992) 51-59.

[12] M. Meyniel, Une condition suffisante d'existence d'un circuit hamiltonien dans un graphe oriente, J. Combin. Theory Ser. B 14 (1973) 137-147.

[13] M. Overbeck-Larisch, Hamiltonian paths in oriented graphs, J. Combin. Theory Ser. B 21 (1) (1976) 76-80.

[14] C. Thomassen, An Ore-type condition implying a digraph to be pancyclic, Discrete Math. 19 (1) (1977) 85-92.

[15] C. Thomassen, Long cycles in digraphs,  Proc. London Math. Soc. (3) 42 (1981) 231-251.

[16] D.R. Woodall, Sufficient conditions for circuits in graphs, Proc. London Math. Soc. 24 (1972) 739-755.

\end{document}